\documentclass[12pt]{amsart}
\usepackage{amsmath, amssymb, amsfonts, amsthm,amscd}
\usepackage{mathrsfs,bbm}
\usepackage{txfonts}
\usepackage[all]{xy}
\usepackage{nicefrac}
\newtheorem{theorem}{Theorem}

\newtheorem{proposition}[theorem]{Proposition}
\newtheorem{conjecture}[theorem]{Conjecture}

\newtheorem{lemma}[theorem]{Lemma}

\theoremstyle{definition}
\newtheorem{definition}[theorem]{Definition}
\newtheorem{remark}[theorem]{Remark}

\newtheorem{example}[theorem]{\bf Example}

\newcommand\SL{\operatorname{SL}}

\renewcommand{\ni}{\noindent}

\title{ Conjectural positivity of Chern-Schwartz-MacPherson classes for Richardson cells}

\author{SHRAWAN KUMAR}
\address{ 
Department of Mathematics, University of North Carolina at Chapel Hill, Chapel Hill, NC 27599-3250, U.S.A.}
\email{shrawan@email.unc.edu}

\date{}

\begin{document}
\maketitle{}

\section{Abstract}  Following some work of  Aluffi-Mihalcea-Sch\"{u}rmann-Su for the CSM classes of Schubert cells and some elaborate computer calculations by R. Rimanyi and L. Mihalcea, I  conjecture that the CSM classes of the Richardson cells expressed in the Schubert basis have nonnegative coefficients. This conjecture was principally motivated by a new product $\square$ coming from the Segre classes in the cohomology of  flag varieties (such that the associated Gr of this product is the standard cup product) and the  conjecture that the structure constants of this new product $\square$ in the standard Schubert basis have alternating sign behavior. I  prove that this conjecture on the sign of the structure constants of $\square$ would follow from my above positivity  conjecture  about the CSM classes of Richardson cells. 

\section{Introduction} Let $G$ be a connected simple algebraic group over $\mathbb{C}$ and let $B$ be a Borel subgroup and $T\subset B$ a maximal torus. Let $W$ be the associated Weyl group with its length function $\ell$. For any $u\in W$, let 
$\mathring{X}_{u}=Bu B/B$ be the Schubert cell and let  $\mathring{X}^{u}:=B^-u B/B$ be the opposite Schubert cell, where 
 $B^-$ is the opposite Borel subgroup of $G$. For $u,v \in W$, let $  \mathring{X}^v_u :=  \mathring{X}_u \cap  \mathring{X}^v$  be the Richardson cell. Let $\left\{\epsilon^{u}\in H^{2 \ell(u)}(G/B, \mathbb{Z}) 
 \right\}_{u\in W }$ be the Schubert cohomology basis of $H^*\left(G/B, \mathbb{Z}\right)$, i.e., $\epsilon^{u}\left(\mu_{X_{v}}\right)=\delta _{u,v}$,  
where $\mu_{X_{v}}$ denotes the fundamental homology class of ${X_{v}}$ in $H_{2\ell(v)}\left(G/B, \mathbb{Z}\right)$. 

 The aim of this paper is to study the  Chern-Schwartz-MacPherson (for short CSM) class $c_{\text{SM}}(\mathring{X}^v_u)$ of Richardson cell  $  \mathring{X}^v_u$ and its connection to some conjecture I learnt from  A. Knutson on a new product $\square$ in the cohomology $H^*(G/B, \mathbb{Z})$ of $G/B$. 
 
 We first prove the following theorem (cf. Theorem \ref{pop5}).
 
 \vskip1ex
\ni 
 {\bf Theorem A.} {\it 
For any $u, v\in  W $, write
$$
c_{{\text{SM}}}\left(  \mathring{X}^v_u\right) = \sum_{w\in W}\, c^{u,v}_ w  \epsilon^{w_ow},\text{ for }c^{u,v}_ w  \in \mathbb{Z},
$$
Then, $c^{u,v}_w = 0$ unless $\ell ( w ) + \ell (u) + \ell (v)$ is even.}
\vskip1ex

We further conjecture the following, which is the main conjecture of this note (cf. Conjecture \ref{conj6}).
 
 \vskip1ex
 \ni
 {\bf Conjecture B.} {\it With the notation as above in Theorem A,
$
c^{u,v}_ w  \geq 0,  \quad \text{for all }u, v,  w .
$}

\vskip1ex

We make the following  conjecture (cf. Conjecture \ref{conj8}) which is a weaker form of the above Conjecture  (as proved in Theorem \ref{thm105}). 
\vskip1ex
\ni
 {\bf Conjecture C.} {\it 
Write, for any $u, v\in W$,
$
c_{\text{SM}} (  \mathring{X}^v_u) = \sum_w d^{u,v}_ w  c_{\text{SM}} (  \mathring{X}_w). 
$
Then, $(-1)^{\ell ( w ) - \ell (u) - \ell(v)} \cdot d^{u,v}_ w  \geq 0$.}
\vskip1ex

Consider the following new product $\square$ in $H^*(G/B, \mathbb{Z})$ coming from the Segre classes, taking the gr of which we recover the standard cup product.
$$
\epsilon^{u}\square \epsilon^{v}=\sum_{\ell ( w )\geq \ell (u)+\ell(v)}\, \chi\left(\mathring{X}\left(u,v, w ;g\right)\right)
\,\epsilon^{ w },
$$
where  $  \mathring{X}\left(u,v, w ;g\right) :=  \mathring{X}_{ w_ou}\cap \mathring{X}^{v} \cap g  \mathring{X}_{ w }$
and $g\in G$ is a general element  such that the intersection  $ {X}\left(u,v, w ;g\right) :=
X_{ w_ou}\cap X^{v} \cap g  X_{ w }$  satisfies the stratified transverse intersection property. We recall the following
conjecture (cf. Conjecture \ref{conj1}).
\vskip1ex
\ni
{\bf Conjecture D.} {\it For any $u, v, w\in W$, 
$$(-1)^{\text{dim}~  \mathring{X}\left(u,v, w ;g\right)}\,
\chi\left( \mathring{X}\left(u,v, w ;g\right)\right)\geq 0,$$
for a general element $g\in G$ such that $ {X}\left(u,v, w ;g\right)$ satisfies
the stratified transverse intersection property.}
\vskip1ex

We prove the following result showing that Conjecture D is equivalent to the Conjecture C (cf. Theorems \ref{thm9} and \ref{thm105} and Proposition \ref{prop2.13}). In fact, we have the following sharper result:
\vskip1ex
\ni
{\bf Proposition E.} For any fixed $u,v\in W$, 
$$
(-1)^{\text{dim}   \mathring{X} (u, v,  w ; g)}
\chi (  \mathring{X} (u, v,  w ; g))  \geq 0 \, \forall w\in W \Longleftrightarrow
\,\,\text{ Conjecture C  is true for $  \mathring{X}^v_{{ w_o}u}$}.$$
Here (as in Conjecture  D) $g\in G$ is a general element such that $ {X}\left(u,v, w ;g\right)$ satisfies
the stratified transverse intersection property.

Thus, the validity of Conjecture B implies Conjecture C.

\vskip3ex
\ni
{\bf Acknowledgements.} I am very grateful to Allen Knutson, Leonardo Mihalcea and Richard Rimanyi for several very helpful conversations/correspondences. I also acknowledge support from the NSF grant DMS-1802328.

 \section{Main Conjecture on CSM classes and Some Results}

We take varieties to be defined over $\mathbb{C}$ and reduced and irreducible.

Let $G$ be a connected simple algebraic group over $\mathbb{C}$ and let $B$ be a Borel subgroup and $T\subset B$ a maximal torus. Let $W$ be the associated Weyl group with its length function $\ell$. For any $u\in W$, define
$$\mathring{X}_{u}=Bu B/B,\quad \mathring{X}^{u}:=B^-u B/B, 
\quad X_{u}=\overline{Bu B/B},\quad X^{u}:=\overline{B^-u B/B},
$$
where $B^-$ is the opposite Borel subgroup of $G$ (i.e., $B^-\cap B= T$). Then, $\mathring{X}_{u}$ (resp. $X_u$) is called the {\it Schubert cell} (resp. {\it Schubert variety}) and  $\mathring{X}^{u}$ (resp. $X^u$) is called the {\it opposite Schubert cell} (resp. {\it opposite Schubert variety}).  

Let $\left\{\epsilon^{u}\right\}_{u\in W }$ be the Schubert cohomology basis of $H^*\left(G/B, \mathbb{Z}\right)$, i.e., $\epsilon^{u}\left(\mu_{X_{v}}\right)=\delta _{u,v}$,  
where $\mu_{X_{v}}$ denotes the fundamental homology class of ${X_{v}}$ in $H_{2\ell(v)}\left(G/B, \mathbb{Z}\right)$. Thus, $\epsilon^u\in H^{2 \ell(u)}(G/B, \mathbb{Z})$. 
\vskip1ex
For any $u, v\in W$,  define the {\it Richardson cell}
$  \mathring{X}^v_u :=  \mathring{X}_u \cap  \mathring{X}^v$  
and let $X^v_u := X_u\cap X^v$ be the Richardson variety. Then, it is nonempty if and only if $v\leq u$. 
\vskip1ex

\begin{definition} (CSM class) ([M], [F1, Example 19.1.7]) Let $X$ be an algebraic variety over the complex numbers $\mathbb{C}$ and let $\mathscr{C}(X)$ be the group of constructible functions on $X$. Thus, the elements of $\mathscr{C}(X)$ are finite sums of the form $\sum_i c_i\mathbb{I}_{X_i}$, where $c_i \in \mathbb{Z}$,  $X_i\subset X$ are locally closed (irreducible) subvarieties and $\mathbb{I}_{X_i}$ is the characteristic function of the set $X_i$. For a proper morphism of varieties $X\to Y$, one defines a $\mathbb{Z}$-linear push-forward 
$$f_*: \mathscr{C}(X) \to \mathscr{C}(Y),\,\,\,\text{where $f_*(\mathbb{I}_W)(y):= \chi(f^{-1}(y)\cap W)$},$$
for a subvariety $W\subset X$ and $y\in Y$, where $\chi$ denotes the topological Euler-Poincar\'e characteristic. Then, MacPherson proved a conjecture of Deligne and Grothendieck stating that there exists a unique natural transformation $c_*: \mathscr{C} \to \bar{H}_*$ such that if $X$ is smooth, then 
$$c_*(\mathbb{I}_X) = c(T_X)\cap \mu_X,$$
where $T_X$ is the tangent bundle of $X$, $c(T_X)$ is the total Chern class of $T_X$, $\mu_X$ is the fundamental class of $X$ and $\bar{H}_*$ denotes  the Borel-Moore homology with integral coefficients. Here the naturality of $c_*$ means that for any proper morphism $f:X\to Y$ between varieties,
$$ f_*(c_*(\phi))= c_*(f_*(\phi))\in \bar{H}_*(Y),\,\,\,\text{for any $\phi \in \mathscr{C}(X)$}.$$

The {\it Chern-Schwartz-MacPherson (for short CSM) class} of any subvariety $W$ of a variety $X$ is defined by
$$c_{\text{SM}}(W) := c_*(\mathbb{I}_W)\in \bar{H}_*(X).$$
By the additivity of $\chi$, for a complete variety $X$, we get
\begin{equation} \label{eqn101} \int c_{\text{SM}}(W)= \chi(W).
\end{equation}
If $X$ is a smooth variety of dimension $d$, then we can identify 
$$\bar{H}_i(X)\simeq H^{2d-i}(X, \mathbb{Z}), \,\,\,\text{(cf. [Fu2, Lemma 2 in $\S$B.2])}.$$
\end{definition}

\begin{definition} \label{defi2.2} (Segre class) 
Define the {\it Segree-SM class} $s_{\text{SM}}$ of any constructible set $Y \subset  G/B$ by 
$$s_{\text{SM}}(Y):= c_{\text{SM}}(Y)\cdot c(T_{G/B})^{-1}.$$
\end{definition}

The following result is due to [AMSS, Corollaries 1.4 and  7.4].
\begin{theorem}\label{lem3}
For $X=G/B$ and $u\in  W$,  write 
$$
c_{\text{SM}}\left(  \mathring{X}_{u}\right)=\sum_{w\in W} a_{u, w } \epsilon^{ w }.
$$
Then,
$a_{u,w}\geq 0, a_{u, w_ou}=1$ and $a_{u,w}=0$ unless $w\geq w_ou$, where $w_o$ is the longest element of the Weyl group $W$.

Moreover,
$$
s_{\text{SM}}\left(  \mathring{X}_{u}\right)=\sum_{w\in W} 
(-1)^{\ell ( w )-\ell ( w_ou)}\quad a_{u, w }\epsilon^ w .\qed
$$
\end{theorem}

\begin{theorem}\label{pop5}
For any $u, v\in  W $, write
$$
c_{{\text{SM}}}\left(  \mathring{X}^v_u\right) = \sum_{w\in W}\, c^{u,v}_ w  [ X_ w ],\text{ for }c^{u,v}_ w  \in \mathbbm{Z},
$$
where $[X_w] =\epsilon^{w_ow}\in H^{2(\ell(w_o)-\ell (w))}(G/B, \mathbb{Z})$
 is the Poincar\'e dual of the fundamental class $\mu_{X_w}$. 

Then, $c^{u,v}_w = 0$ unless $\ell ( w ) + \ell (u) + \ell (v)$ is even.
\end{theorem}

\begin{proof} By [S, Theorem 1.2] and Definition \ref{defi2.2}, 
\begin{eqnarray*}
c_{{\text{SM}}}(  \mathring{X}^v_u) &=& s_{\text{SM}}(  \mathring{X}_u) \cdot c_{\text{SM}}(  \mathring{X}^v)\\
&=&\left(\sum(-1)^{\ell (\theta_1)-\ell (u)}\bar{a}_{u,\theta_1} [ X_{\theta_1}]\right)\cdot \left(\sum \bar{a}_{ w_ov,\theta_2} [ X_{\theta_2}]\right),\,\,
\text{ by Theorem \ref{lem3}}, \quad(\text{where }\bar{a}_{u,\theta} := a_{u, w_o\theta})
\\
&=&(A_+ - A_-) \cdot (B_+ + B_-),\\
\text{where}\quad
A_+ &:=& \sum_{\ell (\theta_1) - \ell (u)  \in 2\mathbbm{Z}} \bar{a}_{u,\theta_1} [ X_{\theta_1}],\,\,\,
A_- := \sum_{\ell (\theta_1) - \ell (u) - 1 \in 2\mathbbm{Z}} \bar{a}_{u,\theta_1} [ X_{\theta_1}]\\
B_+ &:=& \sum_{\ell (\theta_2) - \ell ( w _ov) \in 2\mathbbm{Z}} \bar{a}_{ w_ov, \theta_2} [ X_{\theta_2}],\,\,\,
B_- := \sum_{\ell (\theta_2) - \ell ( w_ov) - 1 \in 2\mathbbm{Z}} \bar{a}_{ w_ov,\theta_2} [ X_{\theta_2}].\\
\end{eqnarray*}
Thus,
\begin{equation}\label{eq4.1}
c_{{\text{SM}}}\left(  \mathring{X}^v_u\right) = (A_+B_+ - A_-B_-) + (A_+B_- - A_-B_+).
\end{equation}
Similarly, $c_{{\text{SM}}}(  \mathring{X}^v_u) = c_{\text{SM}}(  \mathring{X}_u) \cdot s_{\text{SM}}(  \mathring{X}^v)$. Hence, 
\begin{eqnarray}\label{eq4.2}
c_{\text{SM}}(  \mathring{X}^v_u) &=& (A_+ + A_-) \cdot (B_+ - B_-)\nonumber\\
&=& (A_+B_+ - A_-B_-) + (A_-B_+ - A_+B_-).
\end{eqnarray}
Comparing the equations (\ref{eq4.1}) and (\ref{eq4.2}), we get 
$
A_+B_- - A_-B_+ = 0$ and hence
\begin{eqnarray*}
c_{\text{SM}}(  \mathring{X}^v_u) &=& A_+B_+ - A_-B_-\\
&=& \sum_{(-1)^{\ell (\theta_1)} = (-1)^{\ell (u)}\atop (-1)^{\ell (\theta_2)} = (-1)^{\ell ( w_ov)}} \bar{a}_{u,\theta_1} \cdot \bar{a}_{ w_ov,\theta_2} [ X_{\theta_1}] [ X_{\theta_2}]
- \sum_{(-1)^{\ell (\theta_1)} = (-1)^{\ell (u)+1}\atop
(-1)^{\ell (\theta_2)} = (-1)^{\ell ( w_ov)+1}
} \bar{a}_{u,\theta_1} \cdot \bar{a}_{ w_ov,\theta_2} [ X_{\theta_1}] [ X_{\theta_2}]
 \\
&=& \sum_{\ell ( w ) = \ell (u) + \ell (v) (mod ~2)} c^{u,v}_ w  [ X_ w ].\\
&&
\end{eqnarray*}
This proves the theorem. 
\end{proof}
We make the following {\it main conjecture}  of the paper.
\begin{conjecture}\label{conj6}
{\it With the notation as above in Theorem (\ref{pop5}),
$
c^{u,v}_ w  \geq 0,  \quad \text{for all }u, v,  w .
$}

\end{conjecture}

\begin{remark}\label{remark2.7}
{\rm (a) The above conjecture is valid for $\text{SL}_5$ (and hence for $\text{SL}_4, \text{SL}_3, \text{SL}_2$) as the elaborate computer calculation done by R. Rimanyi shows. Also, the above conjecture is valid for all rank$-2$ groups as shown by L. Mihalcea via computer calculation.

(b) L. Mihalcea informed me on August 4, 2022 that Rui Xiong  recently wrote to him mentioning a similar positivity conjecture.}
\end{remark}
Write
$$c_{\text{SM}}(\mathring{X}_u)= \sum_{w\in W}\, d_{u, w} [X(w)].$$
Then, by Theorem \ref{lem3}, $d_{u, w} = 0$ unless $w\leq u$ and $d_{u, u}=1$. Thus, 
$\{c_{\text{SM}}(\mathring{X}_u)\}_{u\in W}$ provides another basis of $\bar{H}_*(G/B, \mathbb{Z})\simeq H^{2 \text{dim} (G/B)-*}\, (G/B, \mathbb{Z})$. 

We make the following  conjecture which is a weaker form of Conjecture \ref{conj6} (as proved in Theorem \ref{thm105}). 
\begin{conjecture} \label{conj8}
{\it Write, for any $u, v\in W$,
$
c_{\text{SM}} (  \mathring{X}^v_u) = \sum_w d^{u,v}_ w  c_{\text{SM}} (  \mathring{X}_w). 
$
Then, $(-1)^{\ell ( w ) - \ell (u) - \ell(v)} \cdot d^{u,v}_ w  \geq 0$.}
\end{conjecture}

Let $\Phi: H^*(G/B) \to H^*(G/B)$ be the ring automorphism such that $\Phi_{|H^{2i}(G/B)} = (-1)^i \text{Id}_{H^{2i}(G/B)}.$ The following lemma was observed by L. Mihalcea.
\begin{lemma} Write
$$\Phi\left(s_{\text{SM}}(\mathring{X}_u^v)\right)= \sum_{w\in W} e^{u,v}_w \epsilon^w.$$
Then, 
$$(-1)^{\ell(w_ou)+\ell(v)}\, e^{u,v}_w\geq 0,\,\,\text{ for all $w\in W$}.$$
\end{lemma}
\begin{proof} By [S, Theorem 1.2],
\begin{align*} 
s_{\text{SM}}(\mathring{X}_u^v) &= s_{\text{SM}}(\mathring{X}_u) \cdot s_{\text{SM}}(\mathring{X}_{w_ov}) \\
&= (-1)^{\ell(w_ou)}\, \Phi\left(c_{\text{SM}}(\mathring{X}_u)\right) \cdot (-1)^{\ell(v)} \Phi\left(c_{\text{SM}}(\mathring{X}_{w_ov}) \right), \,\,\text{by Theorem \ref{lem3}}\\
&= (-1)^{\ell(w_ou)+\ell(v)}\, \Phi\left(c_{\text{SM}}(\mathring{X}_u) \cdot c_{\text{SM}}(\mathring{X}_{w_ov}) \right), \,\,\text{since $\Phi$ is a  ring homomorphism}.
\end{align*}
This proves the lemma by using the first part of Theorem \ref{lem3}.
\end{proof}

\section{A new product in $H^*(G/B)$}

I learnt of the deformed  product $\square$ in $H^*(G/B, \mathbb{Z})$ coming from the Segre classes from A. Knutson (cf. [KZ], [AMSS]). (Taking the gr of $\square$, we clearly recover the standard cup product.)

\begin{definition}  {\rm Define
$$
\epsilon^{u}\square \epsilon^{v}=\sum_{\ell ( w )\geq \ell (u)+\ell(v)}\, \chi\left(\mathring{X}\left(u,v, w ;g\right)\right)
\,\epsilon^{ w },
$$
where we abbreviate $  \mathring{X}\left(u,v, w ;g\right) :=  \mathring{X}_{ w_ou}\cap \mathring{X}^{v} \cap g  \mathring{X}_{ w }$
and $g\in G$ is a general element of $G$ such that the intersection  $ {X}\left(u,v, w ;g\right) :=
X_{ w_ou}\cap X^{v} \cap g  X_{ w }$ is transverse as complex  Whitney stratified  subsets, i.e., the transversality holds at each strata of the triple intersection. We call such an intersection  {\it stratified transverse intersection}. Then, \, $\chi\left( \mathring{X}\left(u,v, w ;g\right)\right)$ does not depend upon the choice of a general element $g\in G$  such that $ {X}\left(u,v, w ;g\right)$ satisfies the stratified transverse intersection property.}
\end{definition}

\begin{example} For $G=\SL_2$, 
	$$\epsilon^e\square \epsilon^e= \epsilon^e-\epsilon^s,\,\,\text{
		where $s$ is the (only) simple reflection}.$$ 
\end{example}

Moreover, I learned of the following non-negativity conjecture from A. Knutson (cf. [KZ, Conjecture after Corollary 1 on Page 43], though they specifically conjectured it only for $4$-step flag manifolds and proved it for $\leq 3$-step flag manifolds; also 2022 ICM talk by A. Knutson). Further, L. Mihalcea informed me that he also has made this conjecture for general $G/B$ and spoke about it in some of the conferences/colloquia.

\begin{conjecture} \label{conj1} {\it For any $u, v, w\in W$, 
$$(-1)^{\text{dim}~  \mathring{X}\left(u,v, w ;g\right)}\,
\chi\left( \mathring{X}\left(u,v, w ;g\right)\right)\geq 0,$$
for a general element $g\in G$ such that $ {X}\left(u,v, w ;g\right)$ satisfies
the stratified transverse intersection property.}
\end{conjecture}

As a special case of Identity \ref{eqn101}, we isolate the following:
\begin{lemma}  \label{lem2.4}
$$
\int_{G/B}c_{\text{SM}}\left(  \mathring{X}\left(u,v, w;g\right)\right)=\chi\left(  \mathring{X}\left(u,v, w; g\right)\right).
$$
\qed
\end{lemma}

The following result is a special case of  [S, Theorem 1.2]. (See an analogous result [O, Proposition 3.8].)
\begin{theorem}\label{thm2}
Assuming the intersection
$  {X}\left(u,v, w ;g\right)$ is stratified transverse, we have 
\begin{eqnarray*}
s_{\text{SM}} \left(\mathring{X}\left(u,v, w;g\right) \right)&=&s_{\text{SM}}\left(  \mathring{X}_{ w_ou}\right)\cdot s_{\text{SM}}\left( \mathring{X}^{v}\right) \cdot s_{\text{SM}}\left(g  \mathring{X}_{ w }\right)\\
&=&s_{\text{SM}}\left(  \mathring{X}_{ w_ou}\right)\cdot s_{\text{SM}}\left( \mathring{X}^{v}\right) \cdot s_{\text{SM}}\left(  \mathring{X}_{ w }\right).\qed
\end{eqnarray*}
\end{theorem}
\begin{theorem} \label{newthm4} For any $u, v, w\in W$ and general $g$, we have
\begin{eqnarray*}
\chi\left(\mathring{X}\left(u,v, w ;g\right)\right)&=&(-1)^{\text{dim}~  \mathring{X}\left(u,v, w ;g\right)}
\,\sum_{\ell (u_1)+\ell (v_1)+\ell ( w _1)=\ell ( w_o)}\, (-1)^{\ell (u)-\ell (u_1)}a_{ w_ou,u_1}a_{ w_ov,v_1}
a_{ w , w _1}
\cdot \,\int_{G/B}\epsilon^{u_1}\cdot \epsilon^{v_1}\cdot \epsilon^{ w _1}.
\end{eqnarray*}
\end{theorem}
\begin{proof}
Let $\left[c_{\text{SM}}(Y)\right]_{\text{top}}$ denote the top degree component of $c_{\text{SM}}(Y)$ for any locally closed subvariety $Y \subset G/B$, i.e., 
writing
$$
c_{\text{SM}}(Y)=\sum_{ w }a_{ w }(Y)\epsilon^{ w }, \quad\quad
\left[c_{\text{SM}}(Y)\right]_{\text{top}}=a_{ w_o }(Y).
$$
There is a similar meaning for $\left[s_{\text{SM}}(Y)\right]_{\text{top}}$. 

Then, by Theorem \ref{thm2} and 
 Definition \ref{defi2.2}, 
\begin{eqnarray*}
c_{\text{SM}}\left(  \mathring{X}\left(u,v, w ;g\right)\right) &=&c_{\text{SM}}\left(  \mathring{X}_{ w_ou}\right)\cdot s_{\text{SM}}\left( \mathring{X}_{ w_ov}\right) \cdot s_{\text{SM}}\left(  \mathring{X}_{ w }\right)\\
&=& \left (\sum_{u_1\geq u} a_{ w_ou,u_1}\epsilon^{u_1}   \right)
\cdot \left(   \sum_{v_1\geq v}(-1)^{\ell (v_1)-\ell (v)} a_{ w_ov,v_1}\epsilon^{v_1} \right)\cdot 
\left(   \sum_{w _1\geq  w_o w }(-1)^{\ell ( w _1)-\ell ( w_o w )} a_{ w , w _1}\epsilon^{ w _1} \right),\\ &&\qquad\qquad\qquad\qquad\qquad\qquad\qquad\qquad\,\,\text{by Theorem \ref{lem3}}\\
&=&\sum_{u_1,v_1, w _1}(-1)^{\ell (v_1)-\ell (v)+\ell ( w _1)-\ell ( w_o w )}
\,\left( a_{ w_ou,u_1}~a_{ w_ov,v_1}~a_{ w , w _1}\epsilon^{u_1}\cdot\epsilon^{v_1}\cdot\epsilon^{ w _1}\right).
\end{eqnarray*}
Thus, by Lemma \ref{lem2.4},
\begin{eqnarray*}
\chi\left(  \mathring{X}\left(u,v, w ;g\right)\right)&=&
\left[c_{\text{SM}} \left(   \mathring{X}(u,v, w ;g)\right)  \right]_{\text{top}}\\
&=&\sum_{u_1,v_1, w _1:\atop
\ell (u_1)+\ell (v_1)+\ell ( w _1)=\ell ( w_o)}(-1)^{\ell ( w )-\ell (u_1)-\ell (v)}
\, a_{ w_ou,u_1}a_{ w_ov,v_1}
a_{ w , w _1}[\epsilon^{u_1}\cdot \epsilon^{v_1}\cdot \epsilon^{ w _1}]_{\text{top}}\\
&=&(-1)^{\text{dim}~  \mathring{X}\left(u,v, w ;g\right)}
\,\sum_{u_1,v_1, w _1:
\atop \ell (u_1)+\ell (v_1)+\ell ( w _1)=\ell ( w_o)}(-1)^{\ell (u) -\ell (u_1)}
a_{ w_ou,u_1}a_{ w_ov,v_1}
a_{ w , w _1}[\epsilon^{u_1}\cdot \epsilon^{v_1}\cdot \epsilon^{ w _1}]_{\text{top}},
\end{eqnarray*}
since 
\setcounter{equation}{0}
\begin{equation} \label{eqn2.11}\text{dim}  \, \mathring{X} (u, v,  w ; g) = \ell ( w ) - \ell (u) - \ell (v).
\end{equation} 
\end{proof}
\section{Conjecture \ref{conj6} implies Conjecture \ref{conj8}}
\begin{theorem}\label{thm9}
Assuming the validity of Conjecture \ref{conj6}, we get for any $u,v, w\in W$, 
$$
(-1)^{\text{dim}  \mathring{X}(u, v,  w ; g)} 
\chi(  \mathring{X}(u, v,  w ; g)) \geq 0,
$$
for any $g$ as in Conjecture  \ref{conj1}. Thus, Conjecture \ref{conj1} is valid assuming Conjecture \ref{conj6}.
\end{theorem}

\begin{proof} By [S, Theorem 1.2] (denoting $\bar{a}_{w, \theta_2}= a_{w, w_o\theta_2}$),
\begin{eqnarray}
c_{\text{SM}} (  \mathring{X} (u, v,  w ; g)) &=& c_{\text{SM}} (  \mathring{X}^v_{ w_ou}) \cdot s_{\text{SM}} (  \mathring{X}_ w )\label{eqn201}\\
 &=& \left(\sum_{\theta_1:\atop \ell (\theta_1) + \ell ( w_ou) + \ell (v)\text{ is even}} c_{\theta_1}^{ w_ou,v} [ X_{\theta_1}]
\right) \cdot 
\, \sum_{\theta_2} (-1)^{\ell (\theta_2)-\ell ( w )} \bar{a}_{w,{\theta_2}} [ X_{\theta_2}],\nonumber \\
&&\qquad\qquad\qquad\qquad \,\,\text{by Theorems \ref{lem3} and \ref{pop5}}
\nonumber\\
&=& \sum_{{\theta_1},{\theta_2}:\atop \ell (\theta_1) + \ell ( w_ou) + \ell (v)\text{ is even}} \,(-1)^{\ell (\theta_2) - \ell ( w )} c_{\theta_1}^{ w_ou,v}
\cdot \,\bar{a}_{w,{\theta_2}} [ X_{\theta_1}] [ X_{\theta_2}].
\end{eqnarray}
Thus,
\begin{eqnarray}\label{eqn4}
 \chi (  \mathring{X} (u, v,  w ; g)) &=&
\int_{G/B} c_{\text{SM}} (  \mathring{X} (u, v,  w ; g)),\,\,\,\text{ by Lemma \ref{lem2.4}}\nonumber\\
&=& \text{coefficient of $[ X_e]$ in} ~	 c_{\text{SM}} (  \mathring{X} (u, v,  w ; g))\nonumber\\
&=& \sum_{{\theta_2} = { w_o}{\theta_1}\atop \ell (\theta_1) + \ell ({ w_o}u) + \ell (v)\text{ is even}} (-1)^{\ell (\theta_2) - \ell ( w )} c_{\theta_1}^{{ w_o}u, v} \bar{a}_{ w ,{\theta_2}}\nonumber\\
&=& (-1)^{\ell (u) + \ell (v) + \ell ( w )} \left(\sum_{{\theta_2} = { w_o}{\theta_1}\atop \ell (\theta_1) + \ell ({ w_o}u) + \ell (v)\text{ is even}} c_{\theta_1}^{{ w_o}u, v} \bar{a}_{ w ,{\theta_2}}\right).
\end{eqnarray}
This proves the theorem by the first part of Therorem \ref{lem3} and the identity \eqref{eqn2.11} assuming the validity of Conjecture \ref{conj6}.
\end{proof}

\begin{proposition}\label{prop2.13} For any fixed $u,v\in W$, 
$$
(-1)^{\text{dim}   \mathring{X} (u, v,  w ; g)}
\chi (  \mathring{X} (u, v,  w ; g))  \geq 0\, \forall w\in W \Longleftrightarrow
\,\,\text{ Conjecture \ref{conj8} is true for $  \mathring{X}^v_{{ w_o}u}$}.$$
Here (as in Conjecture  \ref{conj1}) $g\in G$ is a general element such that $ {X}\left(u,v, w ;g\right)$ satisfies
the stratified transverse intersection property.
\end{proposition}

\begin{proof} By Lemma \ref{lem2.4},
\begin{eqnarray*}
\chi (  \mathring{X} (u, v,  w ; g)) &=& \int_{G/B} c_{\text{SM}} (  \mathring{X} (u, v,  w ; g))\\
&=& \left\langle c_{\text{SM}} (  \mathring{X}^v_{{ w_o}u}), s_{\text{SM}} (  \mathring{X}_ w )\right\rangle,\,\,
\text{by \eqref{eqn201}, where $\langle \alpha, \beta\rangle$ is the coefficient of $[ X_e]$ in $\alpha . \beta$}\\
&=& \left\langle \displaystyle\sum_{\theta} d_{\theta}^{{ w_o}u,{v}} c_{\text{SM}} (  \mathring{X}_{\theta}), s_{\text{SM}} (  \mathring{X}_ w )\right\rangle,\,\,\text{where $d_{\theta}^{{ w_o}u,{v}} $ is as in Conjecture \ref{conj8}}
\\
&=& \displaystyle\sum_{\theta} d_{\theta}^{{ w_o}u, v} \left\langle c_{\text{SM}} (  \mathring{X}_{\theta}), s_{\text{SM}} (  \mathring{X}_ w )\right\rangle\\
&=& d_{{ w_o} w }^{{ w_o}u, v}, \,\,\text{by [AMSS, Corollary 5.9] and Theorem \ref{lem3}.}
\end{eqnarray*}
Recall that, conjecturally,
$\chi (  \mathring{X} (u, v,  w ; g))$ has sign 
$$(-1)^{\text{dim}   \mathring{X} (u, v,  w ; g)}
 =\break (-1)^{\ell ( w ) - \ell (u) - \ell (v)}, \,\,\,\text{by equation (\ref{eqn2.11})}.$$  
Also, by Conjecture \ref{conj8}, $d_{{ w_o} w }^{{ w_o}u, v}$  has sign
$
(-1)^{\ell ({ w_o} w ) - \ell ({ w_o}u ) - \ell (v)} = (-1)^{\ell ( w ) - \ell (u) - \ell (v)}.
$
\end{proof}

Combining Theorem \ref{thm9} and Proposition \ref{prop2.13}, we get the following.

\begin{theorem} \label{thm105} Conjecture \ref{conj8} is  implied by the validity of Conjecture \ref{conj6}.
\end{theorem}


\begin{thebibliography}{300}

\bibitem[AM]{AM} P. Aluffi and  L. Mihalcea: Chern-Schwartz-MacPherson classes for Schubert cells in flag manifolds, {\em Compositio Math.} {\bf 152} (2016), 2603-2625.

\bibitem[AMSS]{AMSS} P. Aluffi, L. Mihalcea, J. Sch\"{u}rmann and C. Su: Shadows of Characteristic cycles, Verma modules, and positivity of Chern-Schwartz-MacPherson clasess of Schubert cells, ArXiv: 1709.08697, v2 (2017).

\bibitem[F1]{F1} W. Fulton, {\em Intersection Theory} (Second Edition), Springer (1998).

\bibitem[F2]{F2} W. Fulton, {\em Young Tableaux}, London Mathematical Society Student Texts {\bf 35}, Cambridge University Press (1997).

\bibitem[KZ]{KZ} A. Knutson and P. Zinn-Justin, Schubert puzzles and integrability II: multiplying motivic Segre classes, ArXiv: 2102.00563, v3
(2021). 


\bibitem[M]{M} R. MacPherson, Chern classes for singular algebraic varieties, {\em Annals of Mathematics} {\bf 100} (1974), 423--432.


 \bibitem[MO]{MO}  D. Maulik and A. Okounkov: Quantum groups and quantum cohomology, ArXiv: 1211.1287 (2012).
 

\bibitem[O]{O}  T. Ohmoto: Singularities of maps and characteristic classes, Advanced Studies in Pure Math. {\bf 68}, 191--265 (2016).

\bibitem[S]{S}  J. Sch\"{u}rmann: Chern classes and transversality for singular spaces, In: ``Singularities in Geometry, Topology, Foliations and Dynamics" (ed. J. Cisneres-Molina et al.), Trends in Mathematics, Birkh\"{a}user, pp. 207--231.
(2017).


\end{thebibliography}
\end{document}